\documentclass[reqno,11pt]{amsart}
\newtheorem{rem}{Remark}

\newtheorem{prop}{Proposition}

\newcommand\eps\varepsilon
\newcommand\ph\varphi
\newcommand\kap\varkappa

\usepackage{enumerate}
\usepackage{graphicx}
\usepackage{float}

\begin{document}

\title[Bounded Solutions to the  System of 2-nd Order  ODE ]
{Bounded Solutions to the  System of 2-nd Order  ODE and the Whitney pendulum}

\author[Oleg Zubelevich]{Oleg Zubelevich\\ \\\tt
 Dept. of Theoretical mechanics,  \\
Mechanics and Mathematics Faculty,\\
M. V. Lomonosov moscow State University\\
Russia, 119899, moscow, Vorob'evy gory, MGU \\
 }
%\email{ozubel@yandex.ru}
%\curraddr{2-nd  Krestovskii Pereulok 12-179, 129110, Moscow, Russia}
\date{}
\thanks{Partially supported by grants
 RFBR  12-01-00441.,  Science Sch.-2964.2014.1}
\subjclass[2000]{  37C60,  34C11}
\keywords{Whitney pendulum, bounded solutions, inverted pendulum, Wazewski method}

\begin{abstract}We propose the existence theorem for bounded solutions to the   system of 2-nd order  ODE. Dynamical applications have been considered.
\end{abstract}

\maketitle
%\tableofcontents
\numberwithin{equation}{section}
\newtheorem{theorem}{Theorem}[section]
\newtheorem{lemma}[theorem]{Lemma}
\newtheorem{definition}{Definition}[section]

\section{Introduction} Courant and Robbins in their book \cite{CR} formulated a problem stated up by  H. Whitney. The problem is as follows.

"Suppose a train travels from station $A$ to station $B$ along a straight section of track. The journey need not be of uniform speed or acceleration. The train may act in any manner, speeding up, slowing down, coming to a halt, or even backing up for a while, before reaching
$B$. But the exact motion of the train is supposed to be known in advance; that is, the function $s=w(t)$ is given, where s is the distance of the train from station $A$, and $t$ is the time, measured from the instant of departure. On the floor of one of the cars a rod is pivoted so that it may move without friction either forward or backward until it touches the floor. If it does touch the floor, we assume that it remains on the floor henceforth; this will be the case if the rod does not bounce. Is it possible to place the rod in such a position that, if it is released at the instant when the train starts and allowed to move solely under the influence of gravity and the motion of the train, it will not fall to the floor during the entire journey from $A$ to $B$?"
\begin{figure}[H]
\centering
\includegraphics[width=45mm, height=40 mm]{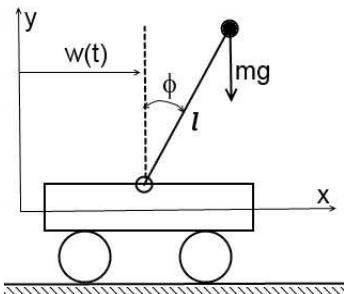}
\caption{The train with the pendulum; $xy$ is the inertial frame. \label{overflow}}
\end{figure}
The authors gave positive answer to this question. Their argument was informal. V. Arnold in \cite{Arn} considered   this problem as  open.

The complete solution to the problem has been given by I. Polekhin in his Ph.D. thesis (unpublished) see also \cite{pal}. He solved the problem by direct application of results from \cite{pol}.

In this article we prove simple and general theorem which implies particularly that there are continuum never-falling solutions to Whitney's problem, we also believe that this theorem describes many other such a type effects.

\section{ Main Theorem}
Introduce several notations. Let $M\subseteq\mathbb{R}^m$ be a domain and $\mathbb{R}_+$ stand for the non-negative reals. A function $$f:H \to\mathbb{R}^m,\quad H=\mathbb{R}_+\times M\times\mathbb{R}^m$$ is $C^2-$smooth in $H$.

The main object of our study is the following system of ordinary differential equations
\begin{equation}\label{fee5}\ddot x=f(t,x,\dot x).\end{equation}Here $x=(x^1,\ldots, x^m)$ is the standard coordinate system in $\mathbb{R}^m$. 

We will use a scalar valued function $F(x),\quad x\in M$ and sets $$D_c=\{x\in M\mid F(x)<c\},\quad \partial D_c=\{x\in M\mid F(x)=c\}.$$

\begin{theorem}\label{dsrger}
Suppose that there exists a  function $F\in C^4(M)$  such that\begin{enumerate}
    \item\label{dfrgd} for some constant $c$ the set $D_c$ is homeomorphic to the open ball of $\mathbb{R}^m$;
\item the set $\overline D_c$ is compact and $\overline D_c\subset M$;
\item\label{fssw4} if $(t,x)\in \mathbb{R}_+\times\partial D_c$ and $\xi\in\mathbb{R}^m$ then the equality $dF(x)[\xi]=0$ implies \begin{equation}\label{swwf4}dF(x)[ f(t,x,\xi)]+ d^2F(x)[\xi,\xi]> 0;\end{equation}
\item\label{sdfsdfddd}  if a solution $x(t),\quad x(0)\in D_c$ to problem (\ref{fee5}) is not defined for all  $t\ge 0$ then it leaves the domain $D_c$ i.e. for some $t'>0$ one has $x(t')\in \partial D_c$.  
\end{enumerate}
Take any   continuous vector field $v(x),\quad x\in M$ such that
$$dF(x)[ v(x)]\Big|_{x\in \partial D_c}\ge 0.$$
Then there exists a point $y\in D_c$ such that system (\ref{fee5}) has a solution\begin{equation}\label{xvfv}x(t)\in C^4(\mathbb{R}_+),\quad x(0)=y,\quad \dot x(0)=v(y)\end{equation} and for all $t\ge 0$ it follows that $x(t)\in D_c$.\end{theorem}
\begin{rem}\label{sdfff}Actually there is no need to demand the domain $D_c$ to be homeomorphic to the ball. The Theorem remains valid if we replace item (\ref{dfrgd}) with the following one: the set $\overline D_c$ is not continuously retractable to its boundary $\partial D_c$. \end{rem}\begin{rem}From item (\ref{fssw4}) it follows that $dF\mid_{\partial D_c}\ne 0$ and thus $\partial D_c$ is a smooth manifold provided $m\ge 2$.\end{rem}

Suppose that  (\ref{fee5}) has the Lagrangian form
$$\frac{d}{dt}\frac{\partial T}{\partial \dot x^i}-\frac{\partial T}{\partial  x^i}=Q_i(t,x,\dot x),\quad T=\frac{1}{2}g_{ij}(x)\dot x^i\dot x^j.$$ 
Here $g_{ij}$ is a Riemann metric in $M$.
Then inequality (\ref{swwf4}) takes the following invariant shape
$$(\nabla_i\nabla_j F(x))\xi^i\xi^j+dF(x)[R(t,x,\xi)]>0,\quad R^i=g^{ij}Q_j .$$

\subsection{Continuum  of Never-falling Solutions to the Whitney Pendulum }By suitable choice of units one can put $g=1,\quad l=1,\quad m=1$. Then the   motion of the pendulum  is described by the equation
\begin{equation}\label{wf4}\ddot\phi=\sin\phi-\ddot w(t)\cos\phi,\quad w\in C^4(\mathbb{R}_+).\end{equation}
Set the initial conditions for this equation as follows
\begin{equation}\label{sfsf}
\phi(0)=\psi,\quad \dot\phi(0)=\lambda((\pi/2)^2-\psi^2).\end{equation}
\begin{prop}
For any $\lambda\in\mathbb{R}$ there exists an angle $ \psi\in (-\pi/2,\pi/2)$ such that problem (\ref{wf4})-(\ref{sfsf}) has a solution $$\phi(t)\in C^4(\mathbb{R}_+),\quad |\phi(t)|<\pi/2\quad\forall t\ge 0.$$\end{prop}
Indeed, this follows immediately from Theorem \ref{dsrger}. One just must put $$F(\psi)=\psi^2,\quad c=(\pi/2)^2,\quad D_c=(-\pi/2,\pi/2),\quad v(\psi)=\lambda((\pi/2)^2-\psi^2)$$ and $\partial D_c=\{\pm\pi/2\}.$

It remains to check  item (\ref{sdfsdfddd}) of the theorem \ref{dsrger}.
But this item follows from the estimate 
$$(\dot\phi(t))^2\le \mathrm{const}+2\Big(\int_0^t\Big|\frac{d^3w}{dt^3}(s)\Big|ds+|\ddot w(t)|\Big).$$
To obtain this formula one must multiply (\ref{wf4}) by $\dot\phi(t)$ and the integrate by parts. The constant  depends on the initial conditions. So that all the solutions to (\ref{wf4}) are defined for all $t\ge 0$.

If the function $\ddot w$ is bounded: $|\ddot w(t)|<C,\quad t\in \mathbb{R}_+$ then the proposition can be made more precise. By $\phi_0\in (0,\pi/2)$ denote the root of the following equation
$$\tan(\phi_0)=C.$$
Then the pendulum has the continuum of solutions $\phi(t)$ such that $$|\phi(t)|<\phi_0\quad \forall t\ge 0.$$ The argument is the same.

\subsection{The ring on the rotating rod}

A   long enough rod rotates around the point $O$ in the vertical plane. The point $O$ is the middle of the rod.
The angle $\phi$ is between the rod and the horizontal axis $x$. The law of rotation $\phi=\phi(t)\in C^3(\mathbb{R}_+)$ is known.
\begin{figure}[H]
\centering
\includegraphics[width=45mm, height=40 mm]{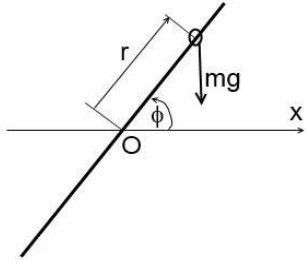}
\caption{The ring on the rotating rod. \label{oow}}
\end{figure}
The rod has a small ring put on it. The ring can slide over the rod without friction. The ring can slide across the point $O$. 

Putting $g=1$ write the equation of motion
$$\ddot r=(\dot\phi(t))^2 r-\sin\big(\phi(t)\big).$$
\begin{prop}
Suppose that  $\dot\phi^2(t)\ge C>0,\quad t\ge 0$. Then one can find an initial position of the ring such that the ring never slides off the rod, provided the rod is sufficiently long.\end{prop}
Indeed, take a function $F(r)=r^2$. Then by the theorem there exists a bounded solution $r(t),\quad |r(t)|<r_*$ and
$r_*$ is chosen  such that
$$r_* C>1.$$

\section{Proof of Theorem \ref{dsrger}}Fix the vector field $v(x)$.  Denote by $x(t,y)$ the solutions with initial conditions (\ref{xvfv}). 

Assume the converse: all the solutions to system (\ref{fee5}) with initial conditions (\ref{xvfv}) leave the domain $D_c$. 

Let  $\tau(y),\quad y\in D_c$ be the time when the solution $x(t,y)$ first time meets $\partial D_c$. That is \begin{equation}\label{4tgrsf}F(x(\tau(y),y))=c.\end{equation} For $y\in\partial D_c$ by definition put $\tau(y)=0$.

\begin{lemma}\label{sdfsf}
If $y\in D_c$ then $dF(x(\tau(y),y))[\dot x(\tau(y),y)]>0$.\end{lemma}
Indeed, assume the converse:  $dF(x(\tau(y),y))[\dot x(\tau(y),y)]\le 0$. Then using the expansion
\begin{align}x(t,y)&=x(\tau(y),y)+\dot x(\tau(y),y)(t-\tau(y))\nonumber\\&+\frac{1}{2}f\big(\tau(y),x(\tau(y),y),\dot x(\tau(y),y)\big)(t-\tau(y))^2+O(t-\tau(y))^3\nonumber
\end{align}
we obtain
\begin{align}F(x(t,y))&=c+dF(x(\tau(y),y))[\dot x(\tau(y),y)]\cdot(t-\tau(y))\nonumber\\&+\frac{1}{2}\Big(d F\big(x(\tau(y),y)\big)[ f\big(\tau(y), x(\tau(y),y),\dot x(\tau(y),y)\big)]\nonumber\\&+
d^2F\big(x(\tau(y),y)\big)[\dot x(\tau(y),y),\dot x(\tau(y),y)]\Big)\cdot(t-\tau(y))^2\nonumber\\&+O(t-\tau(y))^3.\nonumber\end{align}
 By condition (\ref{swwf4}) this formula implies that for small  $\tau(y)-t>0$ it follows that $F(x(t,y))>c$. But this is impossible since for $t<\tau(y)$ the solution $x(t,y)\in D_c$.

The Lemma is proved.
\begin{lemma}\label{sffsdf} The following assertion holds  $\tau(y)\in C(D_c)$.\end{lemma}
Indeed, due to Lemma \ref{sdfsf} this follows from the Implicit function theorem being applied to  equation (\ref{4tgrsf}).

\begin{lemma}\label{sfdd}The following assertion holds $\tau(y)\in C(\overline D_c)$.\end{lemma}
\proof
Fix $\tilde y\in \partial D_c$ i.e. $F(\tilde y)=c$.
We have \begin{equation}\label{sdfs}x(t,y)=y+v(y)t+\frac{1}{2}f(0,y,v(y))t^2+\alpha(t,y),\end{equation}
here $|\alpha(t,y)|\le Kt^3$ if only $y\in D_c$ is close sufficiently to $\tilde y$ and $t$ is small. The constant $K$ is independent on $t,y$.
Substituting formula (\ref{sdfs}) to the equation $F(x(t,y))=c$ we obtain two substantially different situations.

The first one is as follows: $dF(\tilde y)[ v(\tilde y)]>0$ then we have
\begin{equation}\label{dfbd555}\tau(y)=\frac{c-F(y)}{dF( y)[ v(y)]}\Big(1+\lambda(y)\Big),\end{equation}
here $\lambda(y)\to 0$ as $y\to\tilde y$. 

We do not bring detailed proof of formula (\ref{dfbd555}) since  the proof of  analogous fact but just more complicated is  provided below.

From  formula (\ref{dfbd555}) it follows that $\tau(y)\to 0$ as $y\to\tilde y$.

The second case is  $dF(\tilde y)[ v(\tilde y)]=0.$ Due to conditions of the Theorem the last equality implies
$$dF(\tilde y)[ f(0,\tilde y,v(\tilde y))]+ d^2F(\tilde y)[v(\tilde y),v(\tilde y)]> 0.$$
Introduce the notations
\begin{align}
A(y)&=\frac{1}{2}\Big(dF( y)[ f(0, y,v( y))]+ d^2F( y)[v( y),v( y)]\Big),\quad B(y)=dF(y)[v(y)],\nonumber\\
C(y)&=F(y)-c.\nonumber\end{align}
Recall that  we assume  $y\in D_c$ to be close to $\tilde y$. So that $C(y)<0$ and $A(y)\ge c'>0$ with some constant $c'$.

Equation (\ref{4tgrsf}) takes the form
\begin{equation}\label{sdfsdfff}
A(y)\tau^2+B(y)\tau+C(y)+\gamma(\tau,y)=0.\end{equation}
Here $|\gamma (t,y)|\le Kt^3,\quad |\gamma_t(t,y)|\le Kt^2$, constant $K$ is positive, $t\ge 0$ is small enough; and $B(y),C(y)\to 0$ as $y\to\tilde y$.

We seek for solution to (\ref{sdfsdfff}) in the form
$$\tau(y)=u(y)(1+\xi(y)),\quad u=\frac{-B+\sqrt{B^2-4AC}}{2A}.$$
Note that $u(y)>0,\quad u(y)\to 0$ as $y\to\tilde y$.

The function $\xi$ satisfies the equation
\begin{equation}\label{gfff}
\xi+U(y)\xi^2+p(y,\xi)=0,\end{equation}
here \begin{align}\label{egeg}U(y)&=\frac{1}{2}\Big(1-\frac{B(y)}{\sqrt{B^2(y)-4A(y)C(y)}}\Big),\nonumber\\ p(y,\xi)&=\frac{\gamma\big(u(y)(1+\xi),y\big)}{u(y)\sqrt{B^2(y)-4A(y)C(y)}}.\nonumber\end{align} 
It is easy to see that $p(y,\xi)\to 0$ as $y\to\tilde y$ and $0\le U(y)\le 1$.

Consider the equation
$$\eta+s\eta^2+p(y,\eta)=0,\quad 0\le s\le 1.$$ This equation implicitly defines  a function $y\mapsto\eta(\cdot),\quad \eta(s)\in C[0,1]$. Denote this function as $\eta(y,s).$  By the Implicit function Theorem we get this function and $ \eta(y,s)\to 0$ uniformly in $s\in[0,1]$ as $y\to\tilde y$.

The solution to equation (\ref{gfff}) takes the form $\xi(y)=\eta(y,U(y)).$ And consequently $\xi(y)\to 0$ as $y\to\tilde y$.

The Lemma is proved.

Now we can prove the Theorem. By Lemma \ref{sfdd} the mapping $y\mapsto x(\tau(y),y)$ is a continuous retraction of the set $\overline D_c$ to its boundary. It is known from geometry that such a retraction does not exist. 

This contradiction proves the Theorem.

 \end{document}